\documentclass[11pt]{amsart}
\usepackage{amsmath,amssymb}

\usepackage{mathrsfs,bbm,eucal}
\theoremstyle{plain}
\newtheorem{thm}{Theorem}
\theoremstyle{plain}
\newtheorem{lem}{Lemma}
\theoremstyle{plain}
\newtheorem{prop}{Proposition}
\theoremstyle{plain}
\newtheorem{cor}{Corollary}
\theoremstyle{definition}
\newtheorem{defn}{Definition}
\theoremstyle{remark}
\newtheorem{rem}{Remark}

\newcommand{\dirac}{\mbox{$\mathcal{D}\!\!\!\!\!\:/\!\;$}}
\newcommand{\dira}{\mbox{$D\!\!\!\!\!\:/\;\!$}}
\newcommand{\spinor}{\mbox{$S\!\!\!\!\!\:/\;\!$}}
\newcommand{\bundle}[1]{\CMcal{#1}}
\newcommand{\Cl}{\mathit{Cl}}
\newcommand{\R}{\mathbbm{R}}
\newcommand{\C}{\mathbbm{C}}

\newcommand{\T}{\CMcal{T}}

\newcommand{\hyp}[1]{\mbox{$#1\mathrm{H}$}}

\begin{document}
\title[Scalar curvature rigidity of $\hyp{\C }^{2n}$]{Scalar curvature rigidity of almost Hermitian manifolds which are asymptotic to $\hyp{\C}^{2n}$}
\author{Mario Listing}
\address{Department of Mathematics, Stony Brook University, Stony Brook, NY 11794-3651, USA}
\email{listing@math.sunysb.edu}
\thanks{Supported by the German Research Foundation}
\begin{abstract}
We show that an almost Hermitian manifold $(M,g)$ of real dimension $4n$ which is strongly asymptotic to $\hyp{\C }^{2n}$ and satisfies a certain scalar curvature bound must be isometric to the complex hyperbolic space. Assuming K\"ahler instead of almost Hermitian this gives the already known rigidity result by H.~Boualem and M.~Herzlich proved in \emph{Ann.~Scuola Norm.~Sup Pisa (Ser. V)}, vol.~1(2).
\end{abstract}
\maketitle

\section{Introduction}
Scalar curvature rigidity of hyperbolic spaces is a frequently studied problem (cf.~\cite{MinO,AnDa,He3,BoHe,List5,List6}). M.~Herzlich showed in \cite{He3} that a strongly asymptotically complex hyperbolic K\"ahler spin manifold $(M^{2m},g)$ of odd complex dimension $m$ and with scalar curvature $\mathrm{scal}\geq -4m(m+1)$ must be isometric to the complex hyperbolic space $\hyp{\C}^m$. In \cite{BoHe} H.~Boualem and M.~Herzlich gave the corresponding result in the even complex dimensional case, but because of a different representation theory, the spin assumption has to be replaced by another topological condition. In \cite{List6} we generalized the odd complex dimensional case in the way that we only assumed almost Hermitian instead of K\"ahler. In particular, we proved that a complete almost Hermitian spin manifold $(M,g,J)$ of odd complex dimension $m$ which is strongly asymptotically complex hyperbolic and satisfies the scalar curvature bound
\[
\mathrm{scal}\geq -4m(m+1)+6\sqrt{2m}|\nabla J|
\]
must be K\"ahler and isometric to the complex hyperbolic space. In this paper we consider the case of even complex dimension $m$. In order to do so, the spin assumption is replaced by the existence of an appropriate complex line bundle which is the associated line bundle of a chosen spin$^c$ structure. Since an almost Hermitian manifold $(M,g,J)$ is already spin$^c$ (cf.~\cite[App.~D]{LaMi}), there is no need for additional topological assumptions on $M$.

\begin{defn}
\rm $(\hyp{\C}^m,g_0)$ denotes the complex hyperbolic space of complex dimension $m$ and holomorphic sectional curvature $-4$, i.e.~$K\in [-4,-1]$, as well as $B_R(q)\subset M$ is the set of all $p\in M$ with geodesic distance to $q$ less than $R$. Let $(M^{2m},g,J)$ be an almost Hermitian manifold, i.e.~$g$ is a Riemannian metric and $J$ is a $g$--compatible almost complex structure. $(M,g,J)$ is said to be \emph{strongly asymptotically complex hyperbolic} if there is a compact manifold $C\subset M$ and a diffeomorphism $f:E:=M-C\to \hyp{\C}^m-\overline{B_R(0)}$ in such a way that the positive definite gauge transformation $A\in \Gamma (\mathrm{End}(TM_{|E}))$ given by
\[
g(AX,AY) = (f^{*}g_{0})(X,Y)\quad g(AX,Y) = g(X,AY)
\]
satisfies:
\begin{enumerate}
\item $A$ is uniformly bounded.
\item Suppose $r$ is the $f^*g_0$--distance to a fixed point, $\nabla ^0$ is the Levi--Civita connection for $f^*g_0$ and $J_0$ is the complex structure of $\hyp{\C}^m$ pulled back to $E$, then (for some $\epsilon >0$) \[
\left| \nabla ^{0}A\right| +\left| A-\mathrm{Id} \right| +\left| J_0-J\right| \in O(e^{-2(m+1+\epsilon )r}).
\]
\end{enumerate}
\end{defn}

\begin{thm}
\label{rig_thm}
Let $(M,g,J)$ be a complete almost Hermitian manifold of even complex dimension $m=2n$ which is strongly asymptotically complex hyperbolic. Suppose $\frac{\omega }{2\pi }\in \Gamma(\Lambda ^{1,1}M)$ is a closed $2$--form representing the real Chern class of a complex line bundle $\lambda $ which defines a spin$^c$ structure on $M$. If
\begin{equation}
\label{asympt}
2\Omega +\omega \in O(e^{-2(m+1+\epsilon )r}) 
\end{equation}
[$\Omega =g(.,J.)$] and the scalar curvature satisfies
\begin{equation}
\label{ineq}
\mathrm{scal}\geq -4m(m+1)+c_1 |\mathrm{d}^*\Omega |+c_2\Bigl[ |\CMcal{D}^\prime \Omega |+|\CMcal{D}^{\prime \prime }\Omega |\Bigl] +\left| 4\Omega +2\omega \right| ,
\end{equation}
then $(M,g,J)$ is K\"ahler and isometric to $\hyp{\C}^m$.
\end{thm}
In this case $c_1$ and $c_2$ are constants depending on the complex dimension:
\[
c_1:=2\sqrt{\frac{m+1}{m-1}},\qquad c_2:=2\left( m+1-\sqrt{m^2-1}+\frac{2}{\sqrt{m^2-1}}\right) ,
\]
in particular $c_1,c_2\approx 2$ for large $m$. Throughout this paper $\Omega =g(.,J.)$ always denotes the $2$--form associated to $J$, $\mathrm{d}^*$ is formal $L^2$--adjoint of the exterior derivative $\mathrm{d}$ and $\CMcal{D}^\prime +\CMcal{D}^{\prime \prime }$ is the Dolbeault decomposition of $ \CMcal{D}=\mathrm{d}+\mathrm{d}^*$  in $\Lambda ^*(TM)\otimes \C $, i.e.~if $e_1,\ldots ,e_{2m}$ is an orthonormal base, we define $\CMcal{D}^\prime =\sum e_j^{1,0}\cdot \nabla _{e_j}$ and $\CMcal{D}^{\prime \prime }=\sum e_j^{0,1}\cdot \nabla _{e_j}$. Introduce $\CMcal{D}^c:=\mathrm{d}^c+\mathrm{d}^{c,*}$ with $\mathrm{d}^c:=\sum J(e_k)\wedge \nabla _{e_k}$ and $\mathrm{d}^{c,*}:=-\sum J(e_k)\llcorner \nabla _{e_k}$, we obtain $\CMcal{D}^\prime =\frac{1}{2}(\CMcal{D}-\mathbf{i}\CMcal{D}^c)$ as well as $\CMcal{D}^{\prime \prime }=\frac{1}{2}(\CMcal{D}+\mathbf{i}\CMcal{D}^c)$. In particular, we can estimate
\[
|\CMcal{D}^\prime \Omega |+|\CMcal{D}^{\prime \prime }\Omega |\leq |\CMcal{D}\Omega |+|\CMcal{D}^c\Omega |= |\mathrm{d}^*\Omega |+|\mathrm{d}\Omega |+|\mathrm{d}^{c,*}\Omega |+|\mathrm{d}^c\Omega | .
\]
Moreover, using the facts $|\CMcal{D}\eta |\leq \sum _j|\nabla _{e_j}\eta |\leq \sqrt{2m}|\nabla \eta |$ as well as $|\CMcal{D}^c\eta |\leq \sqrt{2m}|\nabla \eta |$ for all $\eta \in \Gamma (\Lambda ^*M)$ we get
\[
|\CMcal{D}^\prime \Omega |+|\CMcal{D}^{\prime \prime }\Omega |\leq 2\sqrt{2m}|\nabla \Omega | .
\]
In complex dimension $m=2$, the constant $c_2$ in inequality (\ref{ineq}) can be improved by setting $c_2:=2(3-\sqrt{3})$. Furthermore, condition (\ref{asympt}) can be replaced (in all dimensions) by one of the following asymptotic assumptions on $g$:\begin{enumerate} \item[(i)] $\mathrm{scal}+4m(m+1)\in O(e^{-2(m+1+\epsilon)r})$\item[(ii)] $|\nabla ^0\nabla ^0A|\in O(e^{-2(m+1+\epsilon )r})$.
\end{enumerate}
Note that (ii) together with the asymptotic assumptions imply (i), and (i) together with inequality (\ref{ineq}) yield (\ref{asympt}). Moreover, assuming that $\Omega $ is a symplectic form (i.e.~$\mathrm{d}\Omega =0$), leads to a much simpler statement of the previous theorem.
\begin{cor}
Let $(M,g,J)$ be a complete almost Hermitian manifold of even complex dimension $m=2n$ which is strongly asymptotically complex hyperbolic, and suppose one of the following conditions:
\begin{enumerate}
\item[(i)] $\mathrm{d}\Omega =0$ and $M$ is diffeomorphic to $\R ^{2m}$.
\item[(i)] $\Omega $ is exact and $M$ is spin.
\item[(ii)] $\mathrm{d}\Omega =0$ and $-\frac{\Omega }{\pi }$ represents the real Chern class of a complex line bundle which defines a spin$^c$ structure on $M$.
\end{enumerate}
If the scalar curvature satisfies
\[
\mathrm{scal}\geq -4m(m+1)+(c_1+2c_2)\sqrt{2m}|\nabla \Omega | ,
\]
then $(M,g,J)$ is K\"ahler and isometric to $\hyp{\C}^{m}$.
\end{cor}
Combining the methods in this paper and in \cite{List6} yield a non--spin version of the result in \cite{List6}.
\begin{prop}
Let $(M,g,J)$ be a complete almost Hermitian manifold of odd complex dimension $m=2n+1$ which is strongly asymptotically complex hyperbolic. Suppose $\lambda $ is the associated complex line bundle of a chosen spin$^c$ structure on $M$, and $\frac{\omega }{2\pi }\in \Gamma(\Lambda ^{1,1}M)$ represents the real Chern class of $\lambda $. If $\omega \in O(e^{-2(m+1+\epsilon )r})$ and the scalar curvature satisfies
\begin{equation*}
\mathrm{scal}\geq -4m(m+1)+2\Bigl[ |\mathrm{d}^*\Omega |+|\CMcal{D}^\prime \Omega |+|\CMcal{D}^{\prime \prime }\Omega |\Bigl] +2\left| \omega \right| ,
\end{equation*}
then $(M,g,J)$ is K\"ahler and isometric to $\hyp{\C}^m$.
\end{prop}

The proof in the complex even dimensional case is much more technical than in the odd dimensional case. We have to modify the ordinary Killing structure in order to show a suitable Bochner--Weitzenb\"ock formula which is necessary to use the non--compact Bochner technique.
\section{Preliminaries}
Let $(M,g,J)$ be an almost Hermitian manifold of complex dimension $m$, then $M$ is a spin$^c$ manifold (cf.~\cite[App.~D]{LaMi}). Suppose $\spinor ^cM$ is a complex spinor bundle of $M$ with associated complex line bundle $\lambda $. We denote by $\gamma $ respectively $\cdot $ the Clifford multiplication on $\spinor ^cM$. $\spinor ^cM$ decomposes orthogonal into
\begin{equation}
\label{dec_spinor}
\spinor ^cM=\spinor ^c_0\oplus \cdots \oplus \spinor ^c_m
\end{equation}
(cf.~\cite{Kir1,LaMi}) where each $\spinor ^c_j$ is an eigenspace of $\gamma (\Omega )$ to the eigenvalue $\mathbf{i}(m-2j)$. Let $\pi _j$ be the orthogonal projections $\spinor ^cM\to \spinor ^c_j$. The decomposition (\ref{dec_spinor}) is parallel (i.e.~$\nabla \pi _j=0 $ for all $j$) if $(g,J)$ is K\"ahler. As usual we introduce $X^{1,0}:=\frac{1}{2}(X-\mathbf{i}J(X))$ and $X^{0,1}:=\frac{1}{2}(X+\mathbf{i}J(X))$. We obtain $\gamma (X^{1,0}):\spinor ^c_j\to \spinor ^c_{j+1}$ as well as $\gamma (X^{0,1}):\spinor ^c_j\to \spinor ^c_{j-1}$, where $\spinor ^c_j=\{ 0\} $ if $j\notin \{ 0,\ldots ,m\}$. If $\nabla ^c$ is a spin$^c$ connection on $\spinor ^cM$, $\dirac ^c$ denotes the corresponding Dirac operator. The $2$--form $\omega $ appearing in the Lichnerowicz formula (cf.~\cite[Thm.~D12]{LaMi}):
\begin{equation}
\label{Lichnerowicz}
(\dirac ^c)^2=\nabla ^*\nabla ^c+\frac{\mathrm{scal}}{4}+\frac{\mathbf{i}}{2}\gamma (\omega )
\end{equation}
will be called the \emph{the curvature $2$--form associated to} $\nabla ^c$ (for the moment $\nabla ^c$ is an arbitrary spin$^c$ connection, but in order to show the main theorem we will consider the canonical spin$^c$ connection induced by the choice of the connection on the complex line bundle $\lambda $). Let $(g,J)$ be K\"ahler and consider the connection $\widehat{\nabla }^0:=\nabla ^c+\frak{A}$ with
\[
\frak{A}_X:= \kappa _1 \gamma (X^{1,0})\pi _{n-1}+\kappa _2\gamma (X^{0,1}) \pi _{n}
\]
$n:=[\frac{m+1}{2}]$. Since $\frak{A}$ is parallel w.r.t.~$\nabla ^c$, we obtain
\begin{equation}
\label{equ_cur}
\begin{split}
\widehat{R}^0_{X,Y}=&\ R^c_{X,Y}+\kappa _1\kappa _2(\gamma (X^{0,1})\gamma (Y^{1,0})-\gamma (Y^{0,1})\gamma (X^{1,0}))\pi _{n-1}\\
&\qquad \quad +\kappa _1\kappa _2(\gamma (X^{1,0})\gamma (Y^{0,1})-\gamma (Y^{1,0})\gamma (X^{0,1}))\pi _n .
\end{split}
\end{equation}
Suppose $\CMcal{R}:\Lambda ^2M\to \Lambda ^2M$ is the Riemannian curvature operator, then the curvature of $\nabla ^c$ satisfies
\begin{equation}
\label{equ_curv}
R^c_{X,Y}=\frac{1}{2}\gamma (\CMcal{R}(X\wedge Y))+\frac{\mathbf{i}}{2}\omega (X,Y) .
\end{equation}
Therefore,
\begin{equation}
\label{equ_xy}
\begin{split}
X^{1,0}\cdot Y^{0,1}-Y^{1,0}\cdot X^{0,1}=\frac{1}{2}\gamma (X\wedge Y+JX\wedge JY-2\mathbf{i}\Omega (X,Y))\\
X^{0,1}\cdot Y^{1,0}-Y^{0,1}\cdot X^{1,0}=\frac{1}{2}\gamma (X\wedge Y+JX\wedge JY+2\mathbf{i}\Omega (X,Y))
\end{split}
\end{equation}
leads to the following proposition:
\begin{prop}
\label{prop1}
Suppose $(M,g,J)$ is a simply connected K\"ahler manifold of constant holomorphic sectional curvature $\kappa $ and complex dimension $m$, then $\widehat{\nabla }^0$ is a flat connection on the subbundle $\bundle{V}=\spinor ^c_{n-1}\oplus \spinor ^c_n\subset \spinor ^cM$, $n=[\frac{m+1}{2}]$, if
\begin{enumerate}\item $m$ is odd, $\kappa =4\kappa _1\kappa _2 $ and $\omega =0$ (i.e.~in particular $M$ is spin)
\item $m$ is even, $\kappa =4\kappa _1\kappa _2$ and $\omega =2\kappa _1\kappa _2\Omega =\frac{1}{2}\kappa \Omega $.
\end{enumerate}
\end{prop}
\begin{proof}
Since $M$ is simply connected, $\widehat{\nabla}^0$ is flat if and only if $\widehat{R}^0=0$. $(M,g,J)$ is of constant holomorphic curvature $\kappa $ if and only if the Riemannian curvature operator satisfies
\[
\CMcal{R}(X\wedge Y)=-\frac{\kappa}{4}(X\wedge Y+JX\wedge JY+2\Omega (X,Y)\Omega ).
\]
Thus, we conclude the claim from $\gamma (\Omega )=\mathbf{i}\sum _j(m-2j)\pi _j$.
\end{proof}

\section{Integrated Bochner--Weitzenb\"ock formula}
In order to show the main theorem, we need a suitable integrated Bochner--Weitzenb\"ock formula. In the case of even complex dimension $m$ it is not possible to show a useful Bochner--Weitzenb\"ock formula for the Killing structure $\frak{A}$ introduced above. However, a minor modification of the Killing structure yields the correct formula. Suppose $(M,g,J)$ is almost Hermitian of even complex dimension $m=2n$, $n\geq 1$. We define $\bundle{V}:=\spinor ^c_{n-1}\oplus \spinor ^c_{n}$, its projection $\mathrm{pr}_\bundle{V}:=\pi _{n-1}+\pi _{n}$ and
\[
\frak{T}_X:=\mathbf{i}\left( \alpha _1X^{1,0}\cdot \pi _{n-1}+\alpha _2X^{0,1}\cdot \pi _{n}+\beta _1X^{1,0}\cdot \pi _{n-2}+\beta _2X^{0,1}\cdot \pi _{n+1} \right)
\]
with (note that $\pi _{n-2}=0$ if $n<2$)
\[
\begin{split}
&\alpha _1=\sqrt{\frac{m-1}{m+1}},\quad  \beta _1=m+1-(m+2)\alpha _1,\\
& \alpha _2=\sqrt{\frac{m+1}{m-1}},\quad   \beta _2=m+1-m\alpha _2,
\end{split}
\]
then $\frak{T}\circ \mathrm{pr}_\bundle{V}$ equals $\frak{A}$ if we set $\kappa _j:=\mathbf{i}\alpha _j$. Moreover, define
\[
\begin{split}
\T &:=-\mathbf{i}(m+2)\alpha _1\pi _{n-1}-\mathbf{i}m\alpha _2\pi _n-\mathbf{i}(m+1)\sum _{j\neq n-1,n}\pi _j\\
&\ = -\mathbf{i}(m+1)+\mathbf{i}\beta _1\pi _{n-1}+\mathbf{i}\beta _2\pi _n .
\end{split} \]
\begin{lem}
$\frak{T}_X+\gamma (X)\T $ is a selfadjoint endomorphism on the complex spinor bundle $\spinor ^cM$ (for every vector field $X$). Moreover, we have (for any orthonormal base $e_1,\ldots ,e_{2m}$)
\begin{equation}
\label{T_prime}
\T ^\prime :=\sum _{j=1}^{2m}(e_j\cdot \frak{T}_{e_j})=\T \circ \mathrm{pr}_\bundle{V}-\mathbf{i}\beta _1(m+4)\pi _{n-2}-\mathbf{i}\beta _2(m+2)\pi _{n+1}.
\end{equation}
\end{lem}
\begin{proof}
Using the facts $(\gamma (X^{1,0}) \pi _j)^*=-\gamma (X^{0,1}) \pi _{j+1}$ and $(\gamma (X^{0,1}) \pi _{j})^*=-\gamma (X^{1,0}) \pi _{j-1}$
we compute
\[
(\frak{T}_X)^*=\mathbf{i}\left( \alpha _1X^{0,1}\pi _n+\alpha _2X^{1,0}\pi _{n-1}+\beta _1X^{0,1}\pi _{n-1}+\beta _2X^{1,0}\pi _n\right)
\]
as well as
\[
\begin{split}
(\gamma (X) \T)^*=-\mathbf{i}(m+1)\gamma (X)+\mathbf{i}\beta _1(X^{0,1}\pi _n+X^{1,0}\pi _{n-2})+\\+\mathbf{i}\beta _2(X^{0,1}\pi _{n+1}+X^{1,0}\pi _{n-1}) .
\end{split}\]
This leads to
\[
\begin{split}
\frak{T}_X+X\cdot \T -&(\frak{T}_X+X\cdot \T )^*=\\&=\mathbf{i}(\alpha _1-\alpha _2+\beta _1-\beta _2)(X^{1,0}\pi _{n-1}-X^{0,1}\pi _n)=0.
\end{split}
\]
The second claim follows from the facts
\begin{equation}
\label{rel1}
\sum _{k=1}^{2m} e_k\cdot e_k^{1,0}=-m+\mathbf{i}\gamma (\Omega )\ \ \text{and}\ \ \sum _{k=1}^{2m} e_k\cdot e_k^{0,1}=-m-\mathbf{i}\gamma (\Omega ).
\end{equation}
\end{proof}

\begin{prop}
Let $(M,g,J)$ be almost Hermitian of even complex dimension $m=2n$. Suppose $\nabla ^c$ is a spin$^c$ connection on $\spinor ^cM$ and $\frak{T}$ as well as $\T $ are given as above. Define the connection $\widehat{\nabla }:=\nabla ^c+\frak{T}$ and the operator $\widetilde{\dirac}:=\dirac ^c +\T$. Then the integrated Bochner--Weitzenb\"ock formula
\[
\int\limits _{\partial N}\left< \widehat{\nabla }_\nu \varphi +\nu \cdot \widetilde{\dirac}\varphi ,\psi \right> = \int\limits _N \left< \widehat{\nabla}\varphi ,\widehat{\nabla }\psi \right> -\left< \widetilde{\dirac}\varphi ,\widetilde{\dirac}\psi \right>+\left< \widehat{\frak{R}}\varphi ,\psi \right>
\]
holds for any compact $N\subset M$ and $\varphi ,\psi \in \Gamma (\spinor ^cM)$. In this case $\nu $ is the outward normal vector field on $\partial N$ and $\widehat{\frak{R}}$ is given by
\[
\begin{split}
\frac{\mathrm{scal}}{4}+\frac{\mathbf{i}}{2}\gamma (\omega )+(m+2)(m-1)\pi _{n-1}+m(m+1)\pi _n+(m+1)^2\mathrm{pr}_{\bundle{V}^\perp}\\-(m+4)\beta _1^2\pi _{n-2}-(m+2)\beta _2^2\pi _{n+1}+\delta \frak{T} +\dira \T
\end{split}\]
while $\omega $ is the curvature $2$--form associated to $\nabla ^c$ [cf.~(\ref{Lichnerowicz})], $\mathrm{pr}_{\bundle{V}^\perp}$ is the projection to the orthogonal complement of $\bundle{V}$ in $\spinor ^cM $, $\delta \frak{T}$ is the divergence of $\frak{T}$, i.e.~$\delta \frak{T}=\sum (\nabla ^c_{e_j}\frak{T})_{e_j}$ and $\dira \T $ is given by $\sum e_j\cdot (\nabla ^c_{e_j}\T )$. Moreover, the boundary operator $\widehat{\nabla }_\nu +\nu \cdot \widetilde{\dirac}$ is selfadjoint.
\end{prop}
\begin{proof}
The selfadjointness of the boundary operator $\widehat{\nabla }_\nu +\nu \cdot \widetilde{\dirac}$ follows immediately from the selfadjointness of $\nabla ^c_\nu +\nu \cdot \dirac ^c$ and $\frak{T}_\nu +\nu \cdot \T $ (since $\nabla ^c$ is a Hermitian connection). The formal $L^2$--adjoint of $\widetilde{\dirac}$ is given by $\widetilde{\dirac}^*=\dirac ^c-\T $. Thus, we can easily verify
\[
\int\limits _N \left< \widetilde{\dirac}\varphi ,\widetilde{\dirac}\psi \right> =-\int\limits _{\partial N}\left< \nu \cdot \widetilde{\dirac}\varphi ,\psi \right> +\int\limits _N\left< \widetilde{\dirac}^*\widetilde{\dirac}\varphi ,\psi \right>
\]
as well as
\[
\begin{split}
\widetilde{\dirac}^*\widetilde{\dirac}=(\dirac ^c)^2+(m+1)^2\mathrm{pr}_{\bundle{V}^\perp}+&\alpha _1^2(m+2)^2\pi _{n-1}+\alpha _2^2m^2\pi _n+\\
&+\dira \T +\sum _{i=1}^{2m}(\gamma (e_i) \T -\T \gamma (e_i))\nabla ^c_{e_i}.
\end{split}\]
Moreover, a straightforward calculation shows
\[ \begin{split}
\int\limits _N\left< \widehat{\nabla }\varphi ,\widehat{\nabla }\psi \right> =&\int\limits _N \left< \nabla ^c\varphi ,\nabla ^c\psi \right> +\left< \nabla ^c\varphi ,\frak{T}\psi \right> +\left< \frak{T}\varphi ,\nabla ^c\psi \right> +\left< \frak{T}\varphi ,\frak{T}\psi \right>\\
=& \int\limits _{\partial N}\left< \nabla ^c_\nu \varphi +\frak{T}_\nu \varphi ,\psi \right> +\int\limits _N \left< \nabla ^*\nabla ^c\varphi ,\psi \right> +\\
&\quad +\int\limits _N\left< \frak{T}\varphi ,\frak{T}\psi \right> -\left< \delta \frak{T}\varphi ,\psi \right> +\Big< \sum _{i=1}^{2m} (\frak{T}_{e_i}^*-\frak{T}_{e_i})\nabla ^c_{e_i}\varphi ,\psi \Big>
\end{split}\]
for all $\varphi ,\psi \in \Gamma (\spinor ^cM)$. Therefore, the Lichnerowicz formula (\ref{Lichnerowicz}) and the fact (previous lemma, $\T ^*=-\T$)
\[
\frak{T}_X^*-\frak{T}_X=\gamma (X) \T - \T \gamma (X)
\]
yields the claim with
\[
\begin{split}
\widehat{\frak{R}}=\frac{1}{4}\mathrm{scal}+\frac{\mathbf{i}}{2}\gamma (\omega )+(m+&1)^2\mathrm{pr}_{\bundle{V}^\perp}+(m+2)^2\alpha _1^2\pi _{n-1}+\\&+m^2\alpha _2^2\pi _n+\dira \T +\delta \frak{T}-\sum _{j=1}^{2m}\frak{T}_{e_j}^*\circ \frak{T}_{e_j}.
\end{split}\]
Use $\pi _j\gamma (X)\pi _{j-1}=\gamma (X^{1,0})\pi _{j-1}$ and $\pi _j\gamma (X)\pi _{j+1}=\gamma (X^{0,1})\pi _{j+1}$ as well as (\ref{rel1}) and $X^{1,0}\cdot X^{1,0}=X^{0,1}\cdot X^{0,1}=0$ to compute
\[
\sum _{j=1}^{2m}\frak{T}_{e_j}^*\circ \frak{T}_{e_j}=(m+2)\alpha ^2_1\pi _{n-1}+m\alpha _2^2\pi _n +(m+4)\beta _1^2\pi _{n-2}+(m+2)\beta _2^2\pi _{n+1}.
\]
Therefore, we obtain $\widehat{\frak{R}}$ from
\[\begin{split}
&(m+2)^2\alpha _1^2-(m+2)\alpha _1^2=(m+2)(m-1)\\
&m^2\alpha _2^2-m\alpha _2^2=m(m+1).
\end{split}\]
\end{proof}
\begin{rem}
\label{rem102}
Since $\widetilde{\dirac}$ is not the Dirac operator of $\widehat{\nabla}$, we made a difference in the notation. However, if $\widehat{\dirac}=\sum  e_j \cdot \widehat{\nabla}_{e_j}$ denotes the Dirac operator of $\widehat{\nabla }$, equation (\ref{T_prime}) yields $\widetilde{\dirac}\circ \mathrm{pr}_\bundle{V}=\widehat{\dirac}\circ \mathrm{pr}_\bundle{V}$. Moreover, if $\varphi $ is a section in $\spinor ^cM$ with $\widehat{\nabla }\varphi =0$ as well as $\widetilde{\dirac }\varphi =0$, then $\varphi $ is a section in $\bundle{V}$ [use the fact $\widehat{\dirac}\varphi =0$ and equation (\ref{T_prime})].
\end{rem}
\begin{lem}
The endomorphism
\[
\begin{split}
-\mathbf{i}\gamma (\Omega )+(m+2)(m&-1)\pi _{n-1}+m(m+1)\pi _n+(m+1)^2\mathrm{pr}_{\bundle{V}^\perp}\\
&-(m+4)\beta _1^2\pi _{n-2}-(m+2)\beta _2^2\pi _{n+1}-m(m+1)\mathrm{Id}
\end{split}\]
is non--negative definite on $\spinor ^cM$ [$\Omega =g(.,J.)$].
\end{lem}
\begin{proof}
Since $\mathbf{i}\gamma (\Omega )=-\sum _j(m-2j)\pi _j$ we
conclude the claim on $\spinor ^c_j$ for $j$ different from $n-2$
and $n+1$. It remains to show that
\[\begin{split}
f_1(m)&:=4+(m+1)-(m+4)\beta _1^2\geq 0\\
f_2(m)&:=-2+(m+1)-(m+2)\beta _2^2\geq 0
\end{split}\]
Both functions are increasing and since $f_1(2)>0$ and $f_2(2)>0$, we get the claim (we only consider the case of complex dimension $m\geq 2$).
\end{proof}
\begin{lem}
Suppose inequality (\ref{ineq}) of the main theorem holds, then at each point of $M$, $\widehat{\frak{R}}$ has no negative eigenvalues: $\widehat{\frak{R}}\geq 0$.
\end{lem}
\begin{proof}
If $\eta $ is a two form, the operator norm of $\gamma (\eta )$ on $\spinor ^cM$ can be estimated by $|\eta |$
\[
|\gamma (\eta )|\leq |\eta |.
\]
Using the last lemma we obtain
\[
\widehat{\frak{R}}\geq \frac{\mathrm{scal}}{4}+m(m+1)-\left| \Omega +\frac{1}{2}\omega \right| -\left| \delta \frak{T}+\dira \T \right| .
\]
Therefore, we have to find an estimate for $\delta \frak{T}+\dira \T $. A straightforward calculation shows
\[
\begin{split}
\delta \frak{T}=&\sum _{j=1}^{2m}(\nabla ^c_{e_j}\frak{T})_{e_j}\\
=&\frac{1}{2}\gamma (\delta J)(\alpha _1\pi _{n-1}-\alpha _2\pi _{n}+\beta _1\pi _{n-2}-\beta _2\pi _{n+1})+\mathbf{i}\sum _{j=1}^{2m}\Big( \alpha _1e_j^{1,0}\nabla ^c_{e_j}\pi _{n-1}\\
&\qquad +\alpha _2e_j^{0,1} \nabla ^c_{e_j}\pi _{n}+\beta _1e_{j}^{1,0} \nabla ^c_{e_j}\pi _{n-2}+\beta _2e_j^{0,1} \nabla ^c_{e_j}\pi _{n+1}\Big) .
\end{split}\]
as well as
\[
\dira \T =\mathbf{i}\sum _{j=1}^{2m}\left( \beta_1e_j\nabla ^c_{e_j}\pi _{n-1}+\beta _2e_j\nabla ^c_{e_j}\pi _n \right) .
\]
Set $\alpha =\alpha _1+\beta _1=\alpha _2+\beta _2$, then $\nabla ^c=\nabla $ on $\Cl _\C(TM)=\mathrm{End}(\spinor ^cM)$ yields
\begin{equation}
\label{estimate}
\begin{split}
|\delta \frak{T}+\dira \T | \leq  \frac{\alpha _2}{2}|\mathrm{d}^*\Omega |&+\alpha \left| \sum \left( e_j^{1,0}\nabla _{e_j}\pi _{n-1}+e_j^{0,1}\nabla _{e_j}\pi _n\right) \right| +\\
&+|\beta _1|\left| \sum \left( e_j^{1,0}\nabla _{e_j}\pi _{n-2}+e_j^{0,1}\nabla _{e_j}\pi _{n-1}\right) \right| +\\
&+|\beta _2|\left| \sum \left( e_j^{1,0}\nabla _{e_j}\pi _{n}+e_j^{0,1}\nabla _{e_j}\pi _{n+1}\right) \right| .
\end{split}
\end{equation}
Thus, we have to estimate $\nabla _X\pi _r$. We conclude from $\pi _{n-k}\gamma (\Omega )=2k\mathbf{i}\pi _{n-k}$
\[
(\nabla _X\pi _{n-k})(2k\mathbf{i}-\gamma (\Omega ))=\pi _{n-k}\gamma (\nabla _X\Omega )
\]
as well as from $\gamma (\Omega )\pi _{n-k}=2k\mathbf{i}\pi _{n-k}$
\[
(2k\mathbf{i}-\gamma (\Omega )) (\nabla _X\pi _{n-k})=\gamma (\nabla _X\Omega )\pi _{n-k} .
\]
Using the facts $\pi _r (\nabla _X\pi _r)\pi _r=0$  for all $r$ and
\[
2k\mathbf{i}-\gamma (\Omega )=\sum _{j\neq {n-k}}c_j\pi _j
\]
with $|c_j|\geq 2$, $|\nabla _X\pi _r|$ can be estimated by $\frac{1}{2}|\nabla _X\Omega |$. Moreover,
\[
\sum _{j=1}^{2m}e_j^{1,0}\cdot (\nabla _{e_j}\pi _{n-k})(2k\mathbf{i}-\gamma (\Omega ))=\pi _{n-k+1}\sum _{j=1}^{2m}\gamma (e_j^{1,0}\cdot \nabla _{e_j}\Omega )
\]
leads to
\[
\Biggl| \sum _{j=1}^{2m}e_j^{1,0}(\nabla _{e_j}\pi _{n-k})\phi \Biggl|\leq \frac{1}{2}\bigl| \gamma (\CMcal{D}^\prime \Omega )\phi \bigl| ,
\]
if $\pi _{n-k}(\phi )=0$, and
\begin{eqnarray*}
\sum _{j=1}^{2m}\gamma (e_j^{1,0}\cdot \nabla _{e_j}\Omega )\pi _{n-k}&=& \sum _{j=1}^{2m}\gamma (e_j^{1,0})(2k\mathbf{i}-\gamma (\Omega ))(\nabla _{e_j}\pi _{n-k})\\
&=& \sum _{j=1}^{2m} ((2k-2)\mathbf{i}-\gamma (\Omega ))\gamma (e_j^{1,0}) (\nabla _{e_j}\pi _{n-k})
\end{eqnarray*}
shows
\[
\Biggl| \sum _{j=1}^{2m}e_j^{1,0}(\nabla _{e_j}\pi _{n-k})\phi \Biggl|\leq \frac{1}{2}\bigl| \gamma (\CMcal{D}^\prime \Omega )\phi \bigl| ,
\]
if $\phi \in \spinor _{n-k}$. In this case we used $\pi _{n-k+1}(e_j^{1,0}\cdot \nabla _{e_j}\pi _{n-k})\pi _{n-k}=0$ and the fact that $(2k-2)\mathbf{i}-\gamma (\Omega )$ has absolute minimal eigenvalue $2$ on $\spinor _{n-k+1}^\perp $. The same procedure yields
\[
\Biggl| \sum _{j=1}^{2m}e_j^{0,1}(\nabla _{e_j}\pi _r)\Biggl|\leq \frac{1}{2}|\CMcal{D}^{\prime \prime }\Omega |.
\]
for all $r$. Thus, we obtain
\[
|\delta \frak{T}+\dira \T |\leq \frac{\alpha _2}{2} |\mathrm{d}^*\Omega |+\frac{\alpha +|\beta _1|+|\beta _2|}{2}\left( |\CMcal{D}^\prime \Omega |+|\CMcal{D}^{\prime \prime }\Omega |\right) ,
\]
and
\[
|\beta _1|+|\beta _2|=\beta _1-\beta _2=\frac{2}{\sqrt{m^2-1}}, \quad \alpha =m+1-\sqrt{m^2-1}
\]
supplies the claim: $\widehat{\frak{R}}\geq 0$. In complex dimension $m=2$ there is a better estimate of $\delta \frak{T}+\dira \T $. Since the decomposition $\spinor ^cM=(\spinor ^cM)^+\oplus (\spinor ^cM)^-$ induced by the volume form is parallel and we have $(\spinor ^cM)^-=\spinor ^c_1$ as well as $(\spinor ^cM)^+=\spinor _0^c\oplus \spinor _2^c$, $\pi _1=\pi _-$ is parallel and we obtain the improvement if $m=2$ from (\ref{estimate}) and the above considerations:
\[
|\delta \frak{T}+\dira \T |\leq \frac{\alpha _2}{2}|\mathrm{d}^*\Omega |+\frac{\alpha }{2}|\CMcal{D}^\prime \Omega |+\frac{|\beta _1|+|\beta _2|}{2}|\CMcal{D}^{\prime \prime }\Omega |.
\]
\end{proof}
\begin{prop}
\label{prop_dirac}
Suppose $(M,g,J)$ is a complete almost Hermitian manifold of complex dimension $m$, $\spinor ^cM$ is a complex spinor bundle of $M$, $\nabla ^c$ is a spin$^c$ connection and $\omega $ is the curvature two form associated to $\nabla ^c$. If the scalar curvature is uniformly bounded with
\begin{equation}
\label{scal_ine}
\frac{\mathrm{scal}}{4}\geq -m(m+1) +\left| \Omega +\frac{1}{2}\omega \right| ,
\end{equation}
then the Dirac operator
\[
\widetilde{\dirac}=\dirac ^c+\T :W^{1,2}(M,\spinor ^cM)\to L^2(M,\spinor ^cM)
\]
is an isomorphism of Hilbert spaces.
\end{prop}
\begin{proof}
First we show that
\[
\overline{\dirac }_\pm :=\dirac ^c\pm \mathbf{i}(m+1):W^{1,2}(M,\spinor ^cM)\to L^2(M,\spinor ^cM)
\]
are isomorphism of Hilbert spaces if the scalar curvature inequality (\ref{scal_ine}) is satisfied. $\overline{\dirac} _\pm $ is bounded on $W^{1,2}$, i.e.~the symmetric bilinear form
\[
B_\pm (\varphi ,\psi ):=\ \int \limits _M \left< \overline{\dirac}_\pm \varphi ,\overline{\dirac}_\pm \psi \right>
\]
is well defined and bounded on $W^{1,2}$. Let $\phi $ be a section in $\spinor ^cM$ with compact support in $M$, then using the Lichnerowicz formula (\ref{Lichnerowicz}) and inequality (\ref{scal_ine}) leads to
\[
\begin{split}
B_\pm (\phi ,\phi )=&\ \int\limits _M \left| \nabla ^c\phi \right| ^2+(m+1)|\phi |^2+ \left< \left( m(m+1)+\frac{\mathrm{scal}}{4}+\frac{\mathbf{i}}{2}\omega \cdot \right)\phi ,\phi \right> \\
\geq &\ \int\limits _M |\nabla ^c\phi | ^2 +(m+1)|\phi |^2-\left< \mathbf{i}\Omega \cdot \phi ,\phi \right> .
\end{split}\]
Therefore, $|\gamma (\Omega ) |\leq m$ on $\spinor ^cM$ implies that $B_\pm $ is coercive, in particular $B_\pm $ is a scalar product on $W^{1,2}$. This proves the injectivity of $\overline{\dirac}_\pm $. The surjectivity of $\overline{\dirac}_\pm $ follows from the Riesz representation theorem and \cite[Thm.~2.8]{GrLa3} (cf.~\cite{AnDa,He3,MinO}). The Dirac operator $\widetilde{\dirac}=\dirac ^c+\T $ is bounded w.r.t.~the $W^{1,2}$--norm, i.e.~$\widetilde{\dirac}$ is well defined and the bilinear form
\[
B(\varphi ,\psi ):=\int \limits _M\left< \widetilde{\dirac }\varphi ,\overline{\dirac}_-\psi \right>
\]
is well defined and bounded on $W^{1,2}(M,\spinor ^cM)$. Using the
definition of $\widetilde{\dirac}$ and the above estimate lead to
\[
\begin{split}
B(\phi ,\phi )=&\ \int\limits _M \left| \overline{\dirac}_- \phi \right| ^2+\mathbf{i}\left< \beta _1\pi _{n-1}(\phi )+\beta _2\pi _n(\phi ),\overline{\dirac}_-\phi \right> \\
\geq & \ \int\limits _M \left( 1-\frac{|\beta _1|+|\beta _2|}{2}\right)  \left| \overline{\dirac}_- \phi \right| ^2- \frac{|\beta _1|}{2}|\pi _{n-1}\phi |^2-\frac{|\beta _2|}{2}|\pi _n\phi |^2 \\
\geq &\ \int\limits _M \left( 1-\frac{|\beta _1|+|\beta _2|}{2}\right) \left( |\nabla ^c\phi |^2+|\phi |^2\right) - \frac{|\beta _1|}{2}|\pi _{n-1}\phi |^2-\frac{|\beta _2|}{2}|\pi _n\phi |^2 .
\end{split}
\]
If $m\geq 2$ one can verify that $|\beta _1|,|\beta _2|< 1$ as
well as $|\beta _1|+\frac{1}{2}|\beta _2|<1$ and $|\beta
_2|+\frac{1}{2}|\beta _1|<1$ (use the fact that $|\beta _1|$ and
$|\beta _1|$ are decreasing and the inequalities hold in case
$m=2$). Thus we conclude that $B$ is coercive on
$W^{1,2}(M,\spinor ^cM)$. In particular, $\widetilde{\dirac}$ has
to be injective and the surjectivity remains to show. Suppose
$\psi \in L^2(M,\spinor ^cM)$, then
\[
l(\phi ):=\int\limits _M\left< \psi ,\overline{\dirac}_-\phi \right>
\]
is a bounded linear functional on $W^{1,2}(M,\spinor ^cM)$. The Lax--Milgram theorem (cf.~\cite[Ch.~5.8]{GilTr}) yields a spinor $\xi \in W^{1,2}(M,\spinor ^cM)$ with
\[
B(\xi ,\phi )=l(\phi )
\]
for all $\phi \in W^{1,2}$. Set $\zeta :=\widetilde{\dirac} \xi -\psi \in L^2$, then $\zeta $ is a weak solution of $\overline{\dirac}_+\zeta $, in this case we used $(\overline{\dirac}_-)^*=\overline{\dirac}_+$. Elliptic theory supplies that $\zeta $ is smooth and since $\dirac ^c\zeta =-\mathbf{i}(m+1)\zeta \in L^2$, theorem 2.8 in \cite{GrLa3} yields $\zeta \in W^{1,2}$. But $\overline{\dirac}_+$ is injective on $W^{1,2}$, i.e.~$\zeta =0$ shows the surjectivity of $\widetilde{\dirac}$.
\end{proof}
\begin{lem}\footnote{Private communication with M.~Herzlich}
\label{lem12}
Suppose $\theta $ is a closed two form on $\hyp{\C}^m$ with $\theta \in O(e^{-\delta r})$, $\delta >3$ ($r$ is the complex hyperbolic distance to a fixed point). Then there is a $1$--form $\eta \in O(e^{-\delta r})$ with $\mathrm{d}\eta =\theta $.
\end{lem}
\begin{proof}
Suppose $X$ is the unit radial vector field of some polar chart on $\hyp{\C}^m-B_{r_0}(0)$, $r_0>0$. Let $\varphi _t$ be the flow of $X$ and define $\eta _1$ by
\[
\eta _1:=-\int\limits _0^\infty X\llcorner \varphi _t^*\theta \ \mathrm{d}t ,
\]
then $\varphi _t^*\in O(e^{3t})$ yields $\mathrm{d}\eta _1=\theta $ as well as $\eta \in O(e^{-\delta r})$ on $\hyp{\C}^m-B_{r_0}(0)$. Suppose $f$ is a cut off function for $B_{r_0}(0)$, i.e.~$f$ is smooth, $f=0$ on $B_{r_0}(0)$ and $f=1$ on $\hyp{\C}^m-B_{r_1}(0)$ for some $r_1>r_0$, then $\theta -\mathrm{d}(f\eta _1)$ is a closed and compact supported two form on $\hyp{\C}^m$. Thus, $\hyp{\C}^m\approx \R ^{2m}$ yields a compact supported $1$--form $\eta _2$ with $\mathrm{d}\eta _2=\theta -\mathrm{d}(f\eta _1)$, and $\eta :=f\eta _1+\eta _2$ satisfies $\mathrm{d}\eta =\theta $ as well as $\eta \in O(e^{-\delta r})$.
\end{proof}
\begin{lem}
\label{lem13}
Suppose $(V,q)$ is a vector space of real dimension $2m$ with a quadratic form $q$ and a $q$--compatible complex structure $J$. Denote by $S=\oplus S_r$ the complex spinor space of $V$ where $S_r$ are induced by the action of the K\"ahler form $\Omega$. Then Clifford multiplication
\[
\gamma _{|\frak{su}}:\Lambda ^{1,1}_0V=\mathfrak{su}(V)\subset \Cl ^c (V,q)\to \mathrm{End}(S_r)
\]
is injective if $0<r<m$ and trivial if $r=0$ or $r=m$. Moreover,
if $\mathrm{Id}$ denotes the identity in $\mathrm{End}(S_r)$,
$\mathbf{i}\cdot \mathrm{Id}$ is not in the image of $\gamma
_{|\frak{su}}$ for all $r=0\ldots m$.
\end{lem}
\begin{proof}
We follow a proof given in \cite{He3}. Clifford multiplication of two forms supplies a representation
\[
\lambda :\Lambda ^{1,1}_{0}(V)\to \mathrm{End}(S_{r}),
\]
since $\theta \cdot S_{r}\subset S_{r}$ for each $\theta \in \Lambda ^{1,1}(V)$. With the inclusion $\frak{su}(m)\subset \frak{so}(2m)$ the representation $\rho :\frak{so}(2m)\to \mathrm{End}(\Cl ^c_{2m})$ (cf.~\cite[Ch.~II (3.1)]{LaMi}) restricts to a representation
\[
\rho :\frak{su}(m)\to \mathrm{End}(\Cl ^c_{2m})=\mathrm{End}(\mathrm{End}(S))
\]
[$\rho $ is the Lie--algebra version of $\mathrm{Ad} :\mathrm{SU}(m)\subset \mathrm{Spin}(2m)\to \mathrm{Aut}(\Cl ^c_{2m})$]. Moreover, the adjoint representation
\[
\mathrm{ad} :\frak{su}(m)\to \mathrm{End}(\Lambda ^{1,1}_{0}(V)).
\]
is irreducible (since $\frak{su}(m)$ is simple). With these definitions $\lambda $ is $\frak{su}(m)$--equivariant, that means
\[
\rho (B)\lambda (\theta _0)=\lambda (\mathrm{ad} (B)\theta _0)
\]
holds on $S_{r}$ for any $B\in \frak{su}(m)$ and $\theta _{0}\in \Lambda ^{1,1}_{0}(V)$. Thus, if $\theta _{0}\in \Lambda ^{1,1}_{0}(V)$ is in the kernel of $\lambda $, $\lambda (\mathrm{ad} (B)\theta _{0})$ vanishes for all $B\in \frak{su}(m)$, in particular $\ker (\lambda )$ is $\frak{su}(m)$--invariant (in the representation $\mathrm{ad}$). In particular the irreducibility of $\mathrm{ad} $ supplies $\ker (\lambda )=\{ 0\} $ or $\ker (\lambda )=\Lambda ^{1,1}_{0}(V)$. The second case appears if and only if $r=0 $ or $r=m$: let $e\in V$ with $q(e)=1$, then Clifford multiplication on $S_{r}$ with $Je\wedge e-\frac{1}{m}\Omega \in \Lambda ^{1,1}_{0}V$ is invertible:
\[
\gamma \Big( e\wedge Je-\frac{1}{m}\Omega \Big) \gamma \Big( Je\wedge e-\frac{1}{m}\Omega \Big) = \mathrm{Id}+\frac{1}{m^2}\gamma (\Omega )^2=\frac{4mr-4r^2}{m^2}\mathrm{Id}_{S_r}.
\]
This implies $\ker (\lambda )=\{ 0\} $ if $0<r<m$. That $\gamma _{|\frak{su}}$ can not be injective on $S_0$ or $S_m$ follows immediately from $\dim S_0=\dim S_m=1$. In order to see the second claim in the case $r\neq 0,m$ we use again the irreducibility of $\mathrm{ad}:\frak{su}(m)\to \mathrm{End}(\frak{su}(m))$. The image of $\gamma _{|\frak{su}(m)}$ is isomorphic (as Lie algebras) to $\frak{su}(m)$. Thus, considering the adjoint transformation
\[
\mathrm{ad}:\mathrm{Im}(\gamma _{|\frak{su}})\to \mathrm{End}(\mathrm{Im}(\gamma _{|\frak{su}}))
\]
yields $\mathrm{ad}(f)\neq 0$ for all $f\in \mathrm{Im}(\gamma _{|\frak{su}})-\{ 0\}$ (otherwise $\mathrm{ad}$ is reducible). In particular, if $\mathbf{i}\cdot \mathrm{Id}$ is contained in $ \mathrm{Im}(\gamma _{|\frak{su}})$, $\mathrm{ad}(\mathbf{i}\cdot \mathrm{Id})=0 $ leads to a contradiction.
\end{proof}
\begin{cor}
\label{cor13}
Clifford multiplication with $\Lambda ^{1,1}M$--forms on $\bundle{V}$ is injective.
\end{cor}
\begin{proof}
Suppose
\[ \eta =\eta _0-\frac{1}{m}\left< \eta ,\Omega \right> \Omega \]
is a $\Lambda ^{1,1}$--form with $\eta _0\in \Lambda ^{1,1}_0$ and $\gamma (\eta )=0$ on $\bundle{V}=\spinor _{n-1}\oplus \spinor _n$. Since $\gamma (\Omega )=0$ on $\spinor _n$, the last lemma yields $\eta _0=0$ ($0<n<m$). Thus, $\gamma (\Omega )=2\mathbf{i}$ on $\spinor _{n-1}$ supplies $\left< \eta ,\Omega \right> =0$ which shows $\eta =0$.
\end{proof}
\section{Proof of the main theorem}

Let $(M,g,J)$ be an almost Hermitian manifold which is strongly asymptotically complex hyperbolic, where $E\subset M$ is supposed to be the Euclidean end of $M$. We denote by $\spinor ^cM$ the considered spin$^c$ bundle as well as by $\nabla ^c$ the canonical spin$^c$ connection determined by the choice of the connection on the complex line bundle $\lambda $ which has curvature $\omega $ (cf.~\cite[Prop.~D11]{LaMi}). The spin connection $\nabla $ is well defined and unique on $\spinor ^cM_{|E}$ and differs from $\nabla ^c$ by an imaginary valued $1$--form: $\mathbf{i}\eta (X):=\nabla ^c_X-\nabla _X$ which satisfies $2\mathrm{d}\eta =\omega $ on $E$. Let $g_0$ be the complex hyperbolic metric on $E$ with K\"ahler structure $\Omega _0:=g_0(.,J_0.)$ [$(g_0,J_0)$ is of constant holomorphic sectional curvature $-4$]. $\nabla ^0$ denotes the Levi--Civita connection for $g_0$ on $TM_{|E}$ as well as the canonical spin connection for $g_0$ on $\spinor ^cM_{|E}$. We conclude from the asymptotic assumptions and condition (\ref{asympt}): $\Omega _0+\frac{1}{2}\omega \in O(e^{-\delta r})$ with $\delta =2(m+1+\epsilon )$. Therefore, lemma \ref{lem12} supplies a $1$--form $\eta _0$ on $E$ (use again a cut off argument) with $\mathrm{d}\eta _0=\Omega _0+\frac{1}{2}\omega $ and $\eta _0\in O(e^{-\delta r})$. Thus, the connection
\[
\nabla ^{0,c}:=\nabla ^0+\mathbf{i}\eta (.)-\mathbf{i}\eta _0(.)
\]
is a spin$^c$ connection on $\spinor ^cM_{|E}$ with associated curvature two form
\[
2\mathrm{d}(\eta -\eta _0)=\omega -2\mathrm{d}\eta _0=-2\Omega _0. \]
Define $\widehat{\nabla }^0:=\nabla ^{0,c}+\frak{A}^0$ on $\spinor ^cM_{|E}$ with $\kappa _j:=\mathbf{i}\alpha _j$ ($j=1,2$), we conclude from proposition \ref{prop1} that $\widehat{\nabla }^0$ is a flat connection on the subbundle $\bundle{V}^0$. We consider the connection $\widehat{\nabla }:=\nabla ^c+\frak{T}$ on $\spinor ^cM$ and show that the restriction of $\widehat{\nabla }$ to $\bundle{V}$ is asymptotic to $\widehat{\nabla }^0$. The gauge transformation $A$ extends to a bundle isomorphism $A:\spinor ^cM_{|E}\rightarrow \spinor ^cM_{|E}$ with (cf.~\cite{AnDa})
\[
\left| \overline{\nabla }\varphi -\nabla ^c\varphi \right| \leq C\left| A^{-1}\right| \left| \nabla ^{0}A\right| \left| \varphi \right| +|\eta _0||\varphi |\, ,\]
where $\overline{\nabla }$ is a connection on $\spinor ^cM_{|E}$ given by $A\nabla ^{0,c}A^{-1}$. Let $\psi _{0}$ be a spinor on $E\subset M$ which is parallel with respect to $\widehat{\nabla }^0$. Set $\psi :=h(A\psi _{0})$ for some cut off function $h$, i.e.~$h=1$ at infinity, $h=0$ in $M-E$ and $\mathrm{supp}(\mathrm{d}h)$ compact. We compute
\[ \begin{split}
\widehat{\nabla }_{X}\psi  = & (Xh)A\psi _{0}+h(\nabla ^c_{X}A\psi _{0}+\frak {T}_{X}(A\psi _{0}))\\
 = & (Xh)A\psi _{0}+h(\nabla ^c_{X}-\overline{\nabla }_{X})A\psi _{0}-hA\frak {A}^0_{X}\psi _{0}+h\frak {T}_{X}A\psi _{0}
\end{split}\]
and thus, the asymptotic assumptions supply
\[
\widehat{\nabla }\psi \in O(e^{(1-\delta)r})\subset L^2(M,T^*M\otimes \spinor ^cM)\]
and
\begin{equation}
\label{41}
\left\langle \widehat{\nabla }_{\nu }\psi +\nu \cdot \widehat{\dirac }\psi ,\psi \right\rangle \in O(e^{(2-\delta )r})\subset L^1(M)
\end{equation}
($|\psi _0|_0^2$ can be estimated by $ce^{2r}$, $|\nu |=1$). Using proposition \ref{prop_dirac} gives a spinor $\xi \in W^{1,2}(M,\spinor ^cM)$ with $\widetilde{\dirac}\xi =\widetilde{\dirac}\psi \in L^2$. In particular $\varphi :=\psi -\xi$ is $\widetilde{\dirac}$--harmonic and non--trivial ($\psi \notin L^2$). Moreover, the selfadjointness of the boundary operator $\widehat{\nabla }_\nu +\nu \cdot \widetilde{\dirac}$ together with (\ref{41}) implies as usual
\[
\liminf _{r\to \infty }\int\limits _{\partial M_r}\left< \widehat{\nabla }_\nu \varphi +\nu \cdot \widetilde{\dirac}\varphi ,\varphi \right> =0
\]
for a non--degenerate exhaustion $\{ M_r\}$ of $M$ (cf.~\cite{AnDa}). Since inequality (\ref{ineq}) gives $\widehat{\frak{R}}\geq 0$, we conclude from the integrated Bochner--Weitzenb\"ock formula:
\[
\int\limits _{\partial M_r}\left< \widehat{\nabla }_\nu \varphi +\nu \cdot \widetilde{\dirac}\varphi ,\varphi \right> \geq \int\limits _{M_r}\left| \widehat{\nabla }\varphi \right| ^2\geq 0 ,
\]
that $\varphi $ is parallel w.r.t.~$\widehat{\nabla }$. Furthermore, $\varphi $ has to be a section of $\bundle{V}=\spinor _{n-1}\oplus \spinor _{n}$ (cf. remark \ref{rem102}), since $0=\widehat{\dirac}\varphi =\dirac ^c\varphi +\T ^\prime \varphi $ and $0=\widetilde{\dirac}\varphi =\dirac ^c\varphi +\T \varphi $. Because $\widehat{\nabla }^0$ is a flat connection in $\bundle{V}^0$, $\bundle{V}$ is trivialized by spinors parallel w.r.t.~$\widehat{\nabla }$. In particular, $\nabla ^c_X$ preserves sections of $\bundle{V}$ which implies that $\mathrm{pr}_\bundle{V}$ is parallel w.r.t.~$\nabla ^c=\nabla $. The complex spinor bundle admits an orthogonal and parallel decomposition $\spinor ^cM=(\spinor ^cM)^+\oplus (\spinor ^cM)^-$ induced from the volume form. If $\pi _+$ as well as $\pi _-$ denote the orthogonal projections of this decomposition, $\pi _{n-1}$ is given by $\mathrm{pr}_\bundle{V}\circ \pi _+$ if $n$ is odd and given by $\mathrm{pr}_\bundle{V}\circ \pi _-$ if $n$ is even. Therefore, $\pi _{n-1}$ as well as $\pi _n=\mathrm{pr}_\bundle{V}-\pi _{n-1}$ have to be parallel and the Killing structure $\frak{A}=\frak{T}\circ \mathrm{pr}_\bundle{V}$ is also parallel w.r.t.~$\nabla ^c$. Moreover, the integrated Bochner--Weitzenb\"ock formula implies $\widehat{\frak{R}}(\varphi )=0$ for all $\varphi \in \Gamma (\bundle{V})$, in particular $\delta \frak{T}=\dira \T =0$ on $\bundle{V}$ supply
\[
0=\frac{\mathrm{scal}}{4}\varphi +\frac{\mathbf{i}}{2}\omega \cdot \varphi +(m+2)(m-1)\pi _{n-1}\varphi +m(m+1)\pi _n\varphi
\]
for all $\varphi \in \bundle{V}$. Thus, lemma \ref{lem13} and $\omega \in \Lambda ^{1,1}M$ yield $\mathrm{scal}=-4m(m+1)$ and $\omega =-2\Omega$. Therefore, $\widehat{R}=0$ on $\bundle{V}$, $\nabla ^c\frak{T}=0$ as well as equations (\ref{equ_cur}), (\ref{equ_curv}) and (\ref{equ_xy}) imply
\begin{equation}
\label{equ_113}
\begin{split}
0=&\ R_{X,Y}^c+[\frak{T}_X,\frak{T}_Y] \\
=&\ \frac{1}{2}\gamma \big (\CMcal{R}(X\wedge Y)+\kappa _1\kappa _2(X\wedge Y+JX\wedge JY+2\Omega (X,Y)\Omega )\big)
\end{split}
\end{equation}
on $\bundle{V}$. From the fact (cf.~\cite{BFGK})
\[
\gamma (\mathrm{Ric}(X)) =2\sum _{i}e_i\cdot R^s_{e_i,X}=\sum _i e_i\cdot \gamma (\CMcal{R}(e_i\wedge X)),
\]
we conclude $\mathrm{Ric}(X)=-2(m+1)X$, i.e.~$g$ is Einstein of scalar curvature $-4m(m+1)$. Inequality (\ref{ineq}) yields $\mathrm{d}^*\Omega =0$ as well as $\CMcal{D}^\prime \Omega =0$ and $\CMcal{D}^{\prime \prime }\Omega =0$. In particular, $\CMcal{D}^\prime +\CMcal{D}^{\prime \prime }=\mathrm{d}+\mathrm{d}^*$ supplies $\mathrm{d}\Omega =0$. Moreover, equation (\ref{equ_113}), $\kappa _1\kappa _2= -1$ and corollary \ref{cor13} show
\begin{equation}
\label{curvature}
\mathrm{pr}_{\Lambda ^{1,1}M}\circ \CMcal{R}(X\wedge Y)=X\wedge Y +JX\wedge JY +2\Omega (X,Y)\Omega .
\end{equation}
Using this equation, the symmetry of the Riemannian curvature tensor and $\Omega \in \Gamma (\Lambda ^{1,1}M)$ lead to
\[
\left<\CMcal{R}(\Omega ),X\wedge Y\right>=\left< \CMcal{R}(X\wedge Y),\Omega \right> =2(m+1)\Omega (X,Y).
\]
Consider the Bochner--Weitzenb\"ock formula on $\Lambda ^2M$:
\[
\triangle =\mathrm{d}^*\mathrm{d}+\mathrm{d}\mathrm{d}^*=\nabla ^*\nabla +\frak{R} ,
\]
then $\frak{R}$ is given by $\mathrm{Ric}+2\CMcal{R}$ (cf.~\cite[Ap.~B]{Se2}), where $\mathrm{Ric}$ acts as derivation on $\Lambda ^2M$. We already know, that $g$ is Einstein, i.e.~$\mathrm{Ric}=-4(m+1)\mathrm{Id} _{\Lambda ^2M}$ supplies $\frak{R}(\Omega )=0$. Moreover, $\mathrm{d}\Omega =0$ and $\mathrm{d}^*\Omega =0$ imply that $\Omega $ is harmonic: $\triangle \Omega =0$, i.e.~we obtain $\nabla ^*\nabla \Omega =0$. Using the fact
\[
0=\triangle |\Omega |^2=\mathrm{d}^*\mathrm{d}|\Omega |^2=2\left< \nabla ^*\nabla \Omega ,\Omega \right> -2\left< \nabla \Omega ,\nabla \Omega \right>
\]
we conclude that $(g,J)$ is K\"ahler. Thus, $\CMcal{R}:\Lambda ^2M\to \Lambda ^{1,1}M$ together with (\ref{curvature}) yield constant holomorphic sectional curvature $-4$ of $(M,g,J)$. Since the end of $M$ is diffeomorphic to $\R ^{2m}-\overline{B_R(0)}$, $M$ must be isometric to $\hyp{\C }^m$.
\bibliographystyle{abbrv}
\bibliography{complex_evem.bbl}

\end{document}